\documentclass{amsart}
\usepackage{amsmath, amssymb, amsthm}
\usepackage{mathtools}
\usepackage{hyperref}
\usepackage{esint}

\numberwithin{equation}{section}

\theoremstyle{definition}
\newtheorem{definition}{Definition}

\newtheorem{theorem}{Theorem}
\newtheorem*{theorem*}{Theorem}
\newtheorem{lemma}{Lemma}
\newtheorem{proposition}{Proposition}
\newtheorem{corollary}{Corollary}

\theoremstyle{remark}
\newtheorem{remark}{Remark}

\newtheorem*{acknowledgement}{Acknowledgement}

\title[Approximation of plurisubharmonic functions]{Approximation of plurisubharmonic functions by logarithms of Gaussian analytic functions}
\author{Kiyoon Eum}
\address{Department of Mathematical Sciences, KAIST, 291 Daehak-ro, Yuseong-gu, Daejeon 34141, South Korea}
\email{kyeum@kaist.ac.kr}

\begin{document}

\begin{abstract}
Let $\Omega$ be a bounded pseudoconvex domain in $\mathbb{C}^N$, and let $u$ be a continuous plurisubharmonic function on $\Omega$. We construct a sequence of Gaussian analytic functions $f_n$ on $\Omega$, associated with $u$, such that $\frac{1}{n}\log|f_n|$ converges to $u$ in $L^1_{loc}(\Omega)$ almost surely, as $n\rightarrow\infty$. Consequently, the normalized zero currents of $f_n$ converge weakly to $dd^c u$. More generally, for every $1\leq k\leq N$, we prove that the normalized currents of simultaneous zeros of $k$ independent copies of $f_n$ converge almost surely to the Bedford--Taylor product $(dd^c u)^k$. We also give a probabilistic proof of the well-known fact that normalized logarithms of the moduli of holomorphic functions are $L^1_{loc}$-dense in the space of plurisubharmonic functions on a pseudoconvex domain.
\end{abstract}

\maketitle

\section{Introduction}
Gaussian random polynomials and Gaussian random analytic functions serve as a means of introducing randomness in complex analysis, allowing for the study of the generic behavior of analytic objects. For instance, the distribution of zeros of random polynomials has been extensively studied in the literature, including \cite{hammersley1956zeros} and \cite{shiffman2003equilibrium}.\par
Another instance of introducing randomness in analysis comes from harmonic analysis. When it is difficult to construct an object satisfying some prescribed properties, it is often easier to define a random object that satisfies these properties almost surely. For related results and methods, see \cite{kahane1985some}.\par
In \cite{bloom2005random} and \cite{bloom2015random}, Bloom and Levenberg considered an asymptotic value of logarithms of random polynomials normalized by their degree and provided a sequence of random polynomials whose normalized logarithms converge to the pluricomplex Green function almost surely. More specifically, their construction involves a sequence of random polynomials of the form
\begin{equation}\label{finite}
    f_n(z)=\sum_{|\nu|\leq n} a_{\nu}^{(n)}p_{\nu}^{(n)}(z)
\end{equation}
where $a_{\nu}^{(n)}$ are i.i.d.\ complex Gaussian random variables and the $p_{\nu}^{(n)}(z)$ form an orthonormal basis of a certain reproducing kernel Hilbert space $\mathcal{H}_n$ of polynomials. Under some condition on $\mathcal{H}_n$ (see (\ref{rk})), they proved that
\begin{equation}\label{bl}
    \left( \limsup_{n\rightarrow\infty}\frac{1}{n}\log |f_n(z)| \right)^*=V_{K}(z) \;\; \text{for all} \; z\in\mathbb{C}^n
\end{equation}
where $V_K(z)$ is the pluricomplex Green function for the regular compact set $K\subset\mathbb{C}^n$ (in the sense of pluripotential theory). Here, $h^*$ means the upper semicontinuous regularization of $h$. They also proved that in (\ref{bl}) the pointwise convergence can be replaced by $L^1_{loc}$ convergence (\cite[Theorem 4.1, Theorem 4.2]{bloom2015random}). See also \cite{bloom2019asymptotic,dauvergne2021necessary} for generalizations of the results in \cite{bloom2005random,bloom2015random}. \par
In this paper, we prove an analogue of (\ref{bl}) for an arbitrary continuous plurisubharmonic (psh for short) function $u$ on $\Omega$ in place of $V_K$. The main ingredient of the construction in \cite{bloom2015random} was the following limiting behavior of the reproducing kernel $B_n(z,w)$ of $\mathcal{H}_n$ (\cite[Proposition 3.1]{bloom2015random}):
\begin{equation}\label{rk}
    \lim_{n\rightarrow\infty}\frac{1}{2n}\log B_n(z,z) = V_K(z) \;\; \text{pointwisely}. 
\end{equation}

For a general psh function $u$ on $\Omega$, Demailly \cite{demailly1992regularization} constructed a sequence of reproducing kernel Hilbert spaces $(\mathcal{H}(nu))_{n\in\mathbb{N}}$ whose reproducing kernels satisfy (\ref{rk}) with $u$ in place of $V_K$. We will define a sequence of Gaussian analytic functions (GAFs) satisfying the analog of (\ref{bl}) using $(\mathcal{H}(nu))_{n\in\mathbb{N}}$, and this will be our main result.

\begin{theorem*}[=Theorem \ref{main}]
Let $u$ be an arbitrary continuous psh function on $\Omega$ and let $(f_n)_{n}$ be a sequence of GAFs defined in Definition \ref{GAF}. Then
\begin{equation*}
(\limsup_{n\rightarrow\infty}\frac{1}{n}\log|f_n(z)|)^*=u(z) \;\; \text{for all}\; z\in\Omega
\end{equation*}
and
\begin{equation*}
\frac{1}{n}\log|f_n(z)|\longrightarrow u  \;\; \text{in} \; L^1_{loc}(\Omega) \;\; \text{as}\; n\rightarrow \infty
\end{equation*}
hold almost surely.
\end{theorem*}

Unlike the case in \cite{bloom2015random}, the property of analytic functions (Theorem \ref{suplp}) is crucially used in the construction of the Gaussian analytic functions and the proof of the results because of the infinite dimensionality of our Hilbert spaces. See the comment before Theorem \ref{main}. 

As a corollary, we give a probabilistic proof of the deterministic theorem of H\"{o}rmander (Corollary \ref{horcor}), stating that the family $\{m^{-1}\log|f| : f\in\mathcal{O}(\Omega),\, m\in\mathbb{N}\}$ is dense in $PSH(\Omega)$ equipped with the $L^1_{loc}(\Omega)$-topology. We denote the space of psh functions on $\Omega$ which are not identically $-\infty$ by $PSH(\Omega)$ and the space of holomorphic functions on $\Omega$ by $\mathcal{O}(\Omega)$. The proof of the main result, together with the submean-value inequality, implies that that the corresponding random variables induced by linear functionals such as $M_v(z,r):=\fint_{B(z,r)}v \,d\lambda$ have exponentially decaying tail probabilities.

Finally, the Poincar\'e--Lelong formula and $L^1_{loc}$-convergence in Theorem \ref{main} imply that the normalized zeros of $f_n$ converge to $dd^c u$ in the sense of currents. We prove that, more generally, for every $1\leq k\leq N$, the normalized simultaneous zero currents of $k$ independent copies of $f_n$ converge almost surely to the Bedford--Taylor product $(dd^c u)^k$ \cite{BT82}.

\begin{theorem*}[=Theorem \ref{zero}]
For every $n\geq1$, let $f_n^1,...,f_n^N$ be independent copies of the GAF in Definition \ref{GAF}. Then for every $1\leq k\leq N$,
\begin{equation}\label{expzero}
\frac{1}{n^k}\mathbb{E} Z_{f_n^1,...,f_n^k}=\left( dd^c\frac{1}{n}\log S_n \right)^k \longrightarrow (dd^cu)^k.
\end{equation}
Moreover,
\begin{equation}\label{justzero}
\frac1{n^k}Z_{f_n^1,...,f_n^k}\longrightarrow (dd^cu)^k
\end{equation}
weakly as currents on $\Omega$, almost surely. 
\end{theorem*}

\begin{remark}[Relation to previous work]
The Gaussian series in Definition \ref{GAF} is the standard Gaussian holomorphic section ensemble associated with an $L^2$-Bergman space. The construction in the possibly infinite dimensional, noncompact setting, as well as the corresponding expected divisor formula, was studied systematically by Drewitz, Liu, and Marinescu \cite[Sections~2 and~3]{Drewitz_Liu_Marinescu_2025}. Their semiclassical equidistribution results concern divisors and are proved under smoothness, uniform positivity, and completeness assumptions. Earlier noncompact equidistribution results were obtained by Dinh, Marinescu, and Schmidt \cite{dinh2012equidistribution}; their section spaces are finite dimensional with controlled growth, and common zeros are discussed under an additional hypothesis in \cite[Remark~3.5]{dinh2012equidistribution}. Related square integrable Gaussian ensembles constructed using Toeplitz operators are studied in \cite{drewitz2024toeplitz}.

For simultaneous zeros, the exact factorization of the expected current is classical in finite dimensional Gaussian ensembles; see, for example, \cite[Proposition~2.1]{bloomshiffman2007zeros} and \cite[Proposition~2.2]{shiffmanzelditch2008variance}.  Almost-sure equidistribution of simultaneous zero currents in finite-dimensional compact or projective settings was proved in \cite[Corollary~1.3]{shiffman2008convergence}; see also \cite{bayraktar2016equidistribution,comanmarinescunguyen2016holder,bojnikgunyuz2026random,guo2026central}. These theorems do not directly give the full statements proved here.  In the present setting the base is a bounded pseudoconvex domain, the weighted Bergman spaces $\mathcal{H}(nu)$ are generally infinite dimensional, and $u$ is only assumed to be continuous and psh, without smoothness or strict plurisubharmonicity; in particular, Bedford--Taylor theory is necessary to define $(dd^cu)^k$.  
\end{remark}

\begin{acknowledgement}
Part of this work was written during the author's visit to the Institute of Mathematics at Jagiellonian University. The author wishes to thank for their hospitality. The author used ChatGPT 5.6 Pro to assist in literature searching and correcting an error in an earlier version of this manuscript. The author independently verified the resulting argument.
\end{acknowledgement}

\section{Probabilistic Preliminaries}
We start with some basic properties of Gaussian random variables. For the definition of complex Gaussian random variables, see \cite{hough2009zeros}, Section 2.1.  We denote by $N(\mu,\sigma^2)$ and $N_{\mathbb{C}}(w,\sigma^2)$ the real and complex Gaussian distributions, respectively.

\begin{proposition}\label{Gau}
Let $a_j \sim N_{\mathbb{C}}(0,\sigma_j^2)$ for each $j\in\mathbb{N}$ and $(a_j)_{j=1}^{\infty}$ are independent. Let $c\in\mathbb{C}$. Then \\
(i) $c a_j\sim N_{\mathbb{C}}(0,|c|^2\sigma_j^2)$ \\
(ii) $\sum_{j=1}^n a_j \sim N_{\mathbb{C}}(0,\sum_{j=1}^n\sigma_j^2)$ \\
(iii) If $\sum_{j=1}^{\infty}\sigma_j^2<\infty$ , then $\sum_{j=1}^{\infty}a_j \sim N_{\mathbb{C}}(0,\sum_{j=1}^{\infty}\sigma_j^2)$
\end{proposition}

As a consequence we have the following lemma.

\begin{lemma}\label{inf}
    Let $(a_j)_{j=1}^{\infty}$ be a sequence of i.i.d.\ standard Gaussian random variables and $(c_j)_{j=1}^{\infty}$ be an $l^2$-sequence of numbers which is not identically zero. Then there exists a constant $C$ independent of $(c_j)_{j=1}^{\infty}$, such that
    \begin{equation*}
\mathbb{P}\left(\frac{|\sum_{j=1}^{\infty} c_j a_j|}{\sqrt{\sum_{j=1}^{\infty}|c_j|^2}}<\frac{1}{n^2}\right) < \frac{C}{n^2}
    \end{equation*}
    for all $n>1$.
\end{lemma}
\begin{proof}
By Proposition \ref{Gau}, we know that 
\begin{equation*}
\frac{\sum_{j=1}^{\infty} c_j a_j}{\sqrt{\sum_{j=1}^{\infty}|c_j|^2}}\sim N_{\mathbb{C}}(0,1).
\end{equation*}
Let $C>0$ be an upper bound for the density function of $N_{\mathbb{C}}(0,1)$ with respect to Lebesgue measure on $\mathbb{C}$. Then 
\begin{equation*}
\mathbb{P}\left(\frac{|\sum_{j=1}^{\infty} c_j a_j|}{\sqrt{\sum_{j=1}^{\infty}|c_j|^2}}<\frac{1}{n^2}\right)=\mathbb{P}\left( |N_{\mathbb{C}}(0,1)|<\frac{1}{n^2}\right)\leq \int_{\{|z|<1/n^2\}} C d\lambda(z) \leq \frac{C\pi}{n^4}
\end{equation*}
where $d\lambda$ denotes Lebesgue measure. The same reasoning also applies to the real Gaussian case. In fact, as the proof shows, $\frac{1}{n}$ instead of $\frac{1}{n^2}$ in the inequality inside of $\mathbb{P}$ would be enough for the complex Gaussian random variables.
\end{proof}

Next, we recall the definition of a Gaussian analytic function. In this paper, we only work on domains in complex Euclidean space. A stochastic process $X(x), x\in\Omega\subset\mathbb{C}^N,$ is a \emph{Gaussian process} if all the finite-dimensional marginals of $X$ are multivariate Gaussian. It is called a \emph{centered Gaussian process} if all of these multivariate Gaussians have mean zero. It is called a \emph{Gaussian analytic function (GAF)} if it is centered and takes values in $\mathcal{O}(\Omega)$. This definition coincides with the one given in \cite{hough2009zeros}, Definition 2.2.1 when $N=1$.

Since a Gaussian random variable is completely determined by its first and second moments, a Gaussian process on $\Omega$ is completely determined by its expectation and covariance functions. In particular, there exists at most one GAF for a given covariance function.

 To prove the main result we need the following lemma. This is a Hilbert space version of Kolmogorov's two-series theorem. For a proof, see \cite[Theorem 2, Chapter 3]{kahane1985some}.
\begin{lemma}\label{kol}
    Let $H$ be a Hilbert space with Hilbert space norm $||\cdot||$. Assume that for each $j\in\mathbb{N}$, $X_j$ is an $H$-valued random vector with $\mathbb{E}(X_j)=0$, $\mathbb{E}||X_j||^2<\infty$. If $(X_j)_{j=1}^{\infty}$ are independent and $\sum_{j=1}^{\infty}\mathbb{E}||X_j||^2<\infty$, then the series  $\sum_{j=1}^{\infty}X_j$ converges in $H$ almost surely. In this case, $\mathbb{E}||\sum_{j=1}^{\infty}X_j||^2=\sum_{j=1}^{\infty}\mathbb{E}||X_j||^2$
\end{lemma}

We end this section by mentioning a few remarks. The first one is that a countable intersection of events of probability one again has probability one. The second one is the Borel--Cantelli lemma, stating that for any sequence of events $\mathcal{A}_n$, if $\sum\mathbb{P}(\mathcal{A}_n)<\infty$, then $\mathbb{P}(\;\mathcal{A}_n \;\text{infinitely often}\;)=0$. Lastly, throughout the paper, almost surely means almost surely in the probability space of the sequence of GAFs or i.i.d. Gaussian random variables, unless otherwise specified. In any case, the meaning will be obvious from the context.

\section{Complex Analytic Preliminaries}

The importance of Lemma \ref{kol} in the study of Gaussian analytic functions lies in the following theorem of complex analysis: 

\begin{theorem}\label{suplp}
    For any compact subset $K$ of $\Omega$ and every relatively compact open neighborhood $\omega$ of $K$ in $\Omega$, there are constants $C_p$ for each $p\in\mathbb{N}$ such that
    \begin{equation*}
    \sup_{K}|f|\leq C_p||f||_{L^{p}(\omega,d\lambda)} 
    \end{equation*}
    for all $f\in\mathcal{O}(\Omega)$.
\end{theorem}
\begin{proof}
    By H\"older's inequality, we only need to prove it for $p=1$. For $p=1$, see \cite[Theorem 2.2.3]{hormander1973introduction}.
\end{proof}

Throughout this paper, we assume $\Omega$ to be a bounded pseudoconvex domain in $\mathbb{C}^N$ and $u$ to be a psh function on $\Omega$. The following estimate of the weighted Bergman kernel will be used. See also the exposition by B\l{}ocki \cite{blocki2009some}.

\begin{proposition}[\protect{\cite[Proposition 3.1]{demailly1992regularization}}]\label{dem}
    Let $u$ be a psh function on $\Omega$. For every $n>0$, let $\mathcal{H}(nu)$ be the Hilbert space of holomorphic functions $f$ on $\Omega$ such that $\int_{\Omega}|f|^2e^{-2nu}d\lambda<\infty$ and let $u_n=\frac{1}{2n}\log \sum_j|\sigma_j^{(n)}|^2$ where $(\sigma_j^{(n)})$ is an orthonormal basis of $\mathcal{H}(nu)$. Then there are constants $L_1, L_2>0$ independent of $n$ such that
    \begin{equation}\label{dineq}
u(z)-\frac{L_1}{n}\leq u_n(z) \leq \sup_{|\zeta-z|<r}u(\zeta)+\frac{1}{n}\log\frac{L_2}{r^N}
    \end{equation}
    for every $z\in\Omega$ and $r<d(z,\partial\Omega)$. In particular, $u_n$ converges to $u$ pointwisely and in $L^1_{loc}(\Omega)$ as $n\rightarrow\infty$.
\end{proposition} 

Note that for each n, $B_n(z,w)=\sum_j\sigma_j^{(n)}(z)\overline{\sigma_j^{(n)}(w)}$ is the reproducing kernels of $\mathcal{H}(nu)$. Since $u$ is locally bounded from above, Theorem \ref{suplp} implies $\sum_j\sigma_j^{(n)}(z)\overline{\sigma_j^{(n)}(w)}$ converges locally uniformly on $\Omega\times\Omega$.

\section{Definition and Convergence Results}
The following lemma provides a convenient way of constructing GAFs when we have explicit information on the covariance function. For example, it applies to Bergman-type kernel functions. This is Lemma 2.2.3 in \cite{hough2009zeros} which is valid for $\Omega\subset\mathbb{C}^N$, $N\geq 1$.

\begin{lemma}\label{pv}
    Let $(\varphi_j)_{j=1}^{\infty}$ be holomorphic functions on $\Omega$ such that $\sum_j |\varphi_j|^2$ converges locally uniformly. Let $(a_j)$ be a sequence of i.i.d. standard complex Gaussian random variables. Then $f(z)=\sum_j a_j\varphi_j(z)$ almost surely converges locally uniformly on $\Omega$ and hence defines a GAF. It has a covariance function $B(z,w)=\sum_j \varphi_j(z)\overline{\varphi_j(w)}$.
\end{lemma}
\begin{proof}
    Let $K$ be a compact subset of $\Omega$ and $\omega$ be a relatively compact open neighborhood of $K$ in $\Omega$. By uniform convergence of $\sum_j |\varphi_j|^2$ on $\omega$, we have $\sum_{j=1}^{\infty}||\varphi_j||_{L^2(\omega)}^2<\infty$. Hence $\sum_j a_j\varphi_j(z)$ satisfies the hypothesis of Lemma \ref{kol} with $a_j\varphi_j$ in place of $X_j$. By Theorem \ref{suplp}, this implies the uniform convergence of $\sum_j a_j\varphi_j(z)$ on $K$. The rest of the statements are immediate.
\end{proof}
Now we can define our GAFs.
\begin{definition}\label{GAF}
    Let $u$ be a continuous psh function on $\Omega$. For each $n\in\mathbb{N}$, let $(\sigma_j^{(n)})_{j=1}^{\infty}$ and $\mathcal{H}(nu)$ be as in Proposition \ref{dem}. We define GAF $f_n$ for each $n\in\mathbb{N}$ by
    \begin{equation}\label{gaf2}
        f_n(z)=\sum_{j=1}^{\infty} a_j^{(n)}\sigma_j^{(n)}(z)
    \end{equation}
    where $(a_j^{(n)})$ is a sequence of i.i.d. standard complex Gaussian random variables. 
    
    Since the invariant sum $\sum_{j=1}^{\infty} \sigma_j^{(n)}(z)\overline{\sigma_j^{(n)}(w)}$ is such an important object in several complex variables, namely a (weighted) Bergman kernel, we would like to have a similar invariance for the sum in (\ref{gaf2}). Indeed, by the remark after the definition of a GAF, this sum is independent of the choice of the orthonormal basis $(\sigma_j^{(n)})_{j=1}^{\infty}$. We denote their covariance functions by 
    \begin{equation*}
        B_n(z,w)=\sum_{j=1}^{\infty} \sigma_j^{(n)}(z)\overline{\sigma_j^{(n)}(w)}
    \end{equation*}
    and use the notation $S_n(z):=\sqrt{B_n(z,z)}=\sqrt{\sum_{j=1}^{\infty}|\sigma_j^{(n)}(z)|^2}$.
\end{definition}
With this definition, we prove our main theorem. Note that except for \emph{Step 1} and \emph{Step 2}, we closely follow the proof of Theorem 4.1 in \cite{bloom2015random}. The main technical difference comes from the fact that in \cite{bloom2015random}, the sum defining (\ref{finite}) is a finite sum, and hence Cauchy--Schwartz inequality can be applied to get an upper bound. However, in our case, the sum has infinitely many summands and $(a_j^{(n)})_{j=1}^{\infty}$ is not $l^2$ almost surely. Hence, we use a different method to get an upper bound, which is inspired by the proof of Lemma \ref{pv}. Namely, we consider $f_n$ as an $L^2(B(z))$-valued random vector for some ball $B(z)$ compactly contained in $\Omega$. The last part of Lemma \ref{kol} provides an explicit formula for the variance of this random vector, which is readily controlled by $S_n$, which in turn is controlled by Proposition \ref{dem}. In this regard, the independence assumption on the coefficients is crucial.

\begin{theorem}\label{main}
    Let $u$ be an arbitrary continuous psh function on $\Omega$ and let $(f_n)_{n}$ be a sequence of GAFs defined in Definition \ref{GAF}. Then
\begin{equation}\label{claim1}
(\limsup_{n\rightarrow\infty}\frac{1}{n}\log|f_n(z)|)^*=u(z) \;\; \text{for all} \; z\in\Omega
\end{equation}
and
\begin{equation}\label{claim2}
\frac{1}{n}\log|f_n(z)|\longrightarrow u  \;\; in \; L^1_{loc}(\Omega) \;\; \text{as}\; n\rightarrow \infty
\end{equation}
hold almost surely.
\end{theorem}
\begin{proof}
\emph{Step 1.---} First we prove that for any $\varepsilon>0$ and $z\in\Omega$, there exists $0<r<d(z,\partial\Omega)$ such that 
\begin{equation*}
    \frac{|f_n(\zeta)|}{S_n(\zeta)}\leq e^{\varepsilon n} 
\end{equation*}
for all but finitely many $n$ on $B(z,r)=\{|\zeta-z|<r\}$ almost surely. By the Borel--Cantelli lemma, it is enough to prove 
\begin{equation}\label{BC}
\sum_{n=1}^{\infty}\mathbb{P}\left(\sup_{B(z,r)}\frac{|f_n(\zeta)|}{S_n(\zeta)}> e^{\varepsilon n} \right)<\infty
\end{equation}
for some $r>0$. Using continuity of $u$, let $0<r_2<d(z,\partial\Omega)$ so that
\begin{equation*}
\text{osc}_{B(z,r_2)}u=\sup_{B(z,r_2)}u-\inf_{B(z,r_2)}u<\varepsilon.
\end{equation*}
Then take $r, r_1$ with $0<r<r_1<r_2$.
By Theorem \ref{suplp}, we have
\begin{equation*}
    \sup_{B(z,r)}\frac{|f_n|}{S_n}\leq \frac{\sup_{B(z,r)}|f_n|}{\inf_{B(z,r)} S_n}\leq \frac{C_2||f_n||_{L^2(B(z,r_1))}}{\inf_{B(z,r)}S_n}.
\end{equation*}
By Chebyshev's inequality and Lemma \ref{kol}, 
\begin{equation*}
\begin{gathered}
    \mathbb{P}\left(\frac{||f_n||_{L^2(B(z,r_1)}}{\inf_{B(z,r)}S_n}>e^{\varepsilon n}\right) \leq \frac{\mathbb{E}||f_n||^2_{L^2(B(z,r_1))}}{e^{2\varepsilon n}(\inf_{B(z,r)}S_n)^2}=\frac{\sum_j \mathbb{E}||a_j^{(n)}\sigma_j^{(n)}||^2_{L^2(B(z,r_1))}}{e^{2\varepsilon n}(\inf_{B(z,r)}S_n)^2} \\
    = \frac{\sum_j ||\sigma_j^{(n)}||^2_{L^2(B(z,r_1))}}{e^{2\varepsilon n}(\inf_{B(z,r)}S_n)^2} \leq \frac{\lambda(B(z,r_1))\sup_{B(z,r_1)}\sum_j|\sigma_j^{(n)}|^2}{e^{2\varepsilon n}(\inf_{B(z,r)}S_n)^2}\\
=\frac{\lambda(B(z,r_1))\sup_{B(z,r_1)}S_n^2}{e^{2\varepsilon n}(\inf_{B(z,r)}S_n)^2}.
\end{gathered}
\end{equation*}

Inequality (\ref{dineq}) implies
\begin{equation}\label{nonzero}
e^{2n\inf_{B(z,r_1)}u-L'_1}\leq e^{2nu(\zeta)-L'_1}\leq S_n^2(\zeta) \leq e^{2n\sup_{B(z,r_2)}u+L'_2}
\end{equation}
on $B(z,r_1)$ for some constants $L'_1, L'_2$ independent of $\zeta$ and $n$. Thus we have
\begin{equation}\label{sum}
\frac{\sup_{B(z,r_1)}S_n^2}{e^{2\varepsilon n}(\inf_{B(z,r)}S_n)^2}\leq e^{2n(\sup_{B(z,r_2)}u-\inf_{B(z,r_1)}u-\varepsilon)+L}\leq e^{-ln+L}
\end{equation}
for some constants $l>0$ and $L$. Since $\sum_n e^{-ln}<\infty$ for $l>0$, this proves (\ref{BC}) and hence the claim. Note we have shown that
\begin{equation}\label{step1}
\mathbb{P}\left(  \sup_{B(z,r)} \frac{|f_n|}{S_n} >e^{\varepsilon n} \right) \leq C e^{-D n} 
\end{equation}
for some $C,D>0$. This estimate will be used in Proposition \ref{subexpthm}.
\par

\emph{Step 2.---} Let $(K_m)_{m=1}^{\infty}$ be a compact exhaustion of $\Omega$ with $\Omega=\bigcup_m K_m$. Fix $K_m$ and $\varepsilon>0$. For each $z\in K_m$, there exist $0<r_z<d(z,\partial\Omega)$ and a random number $N_z<\infty$ such that $\frac{|f_n(\zeta)|}{S_n(\zeta)}\leq e^{\varepsilon n}$ on $B(z,r_z)$ for $n>N_z$ almost surely, by \emph{Step 1}. Since $K_m$ is compact, there exist finitely many points $z_1,...z_M$ in $K_m$ such that 
$K_m\subset\bigcup_i B(z_i,r_{z_i})$. Then for $n>\max_i {N_{z_i}}$, $\frac{|f_n(\zeta)|}{S_n(\zeta)}\leq e^{\varepsilon n}$ on $K_m$ almost surely, and hence
\begin{equation*}
    \frac{1}{n}\log |f_n| \leq \frac{1}{n}\log S_n + \varepsilon
\end{equation*}
on $K_m$ for $n>\max_i {N_{z_i}}$  almost surely. Combined with inequality (\ref{dineq}), this implies $(\frac{1}{n}\log |f_n|)_n$ is almost surely locally uniformly bounded from above. Also by taking the $\limsup$, we get 
\begin{equation*}
\limsup_{n\rightarrow\infty}\frac{1}{n}\log |f_n|\leq u +\varepsilon
\end{equation*}
on $K_m$ almost surely. By letting $\varepsilon=q^{-1}, q=1,2,...$ and taking the countable intersection of the corresponding events, one for each $q$, we get 
\begin{equation}\label{limsup}
\limsup_{n\rightarrow\infty}\frac{1}{n}\log |f_n|\leq u
\end{equation}
on $K_m$ almost surely. Again, taking the countable intersection of the corresponding events, one for each $K_m$, we conclude that the inequality (\ref{limsup}) holds on $\Omega$ almost surely.\par

\emph{Step 3.---} Let $(z_i)_{i=1}^{\infty}$ be a countable dense subset of $\Omega$. For each $z_i$, Lemma \ref{inf} implies
\begin{equation*}
\mathbb{P}\left(\frac{|f_n(z_i)|}{S_n(z_i)}<\frac{1}{n^2}\right)=\mathbb{P}\left(\frac{|\sum_{j=1}^{\infty} a_j^{(n)}\sigma_j^{(n)}(z_i)|}{\sqrt{\sum_{j=1}^{\infty}|\sigma_j^{(n)}(z_i)|^2}}<\frac{1}{n^2}\right) < \frac{C}{n^2}.
\end{equation*}
Note that $S_n(z_i)>0$ by (\ref{nonzero}).
Since $\sum_n 1/n^2 <\infty$, by the Borel--Cantelli lemma,
\begin{equation*}
\frac{|f_n(z_i)|}{S_n(z_i)}\geq\frac{1}{n^2}
\end{equation*}
for all but finitely many $n$ almost surely. Thus we get
\begin{equation}\label{liminf}
\liminf_{n\rightarrow\infty}\frac{1}{n}\log |f_n(z_i)|\geq u(z_i)
\end{equation}
almost surely. Taking the countable intersection of the corresponding events, one for each $z_i$, we prove that (\ref{liminf}) holds for all $z_i$ almost surely. \par

\emph{Step 4.---} Here we prove (\ref{claim1}). In what follows every statement is understood to hold almost surely. Define 
\begin{equation*}
    F(z):=\left( \limsup_{n\rightarrow\infty}\frac{1}{n}\log |f_n(z)| \right)^*.
\end{equation*}
By (\ref{limsup}), 
\begin{equation*}
F\leq u^*=u \;\, on\;\Omega.
\end{equation*}
Also for each $z_i$, 
\begin{equation*}
\limsup_{n\rightarrow\infty}\frac{1}{n}\log |f_n(z_i)|\geq u(z_i)
\end{equation*}
by (\ref{liminf}). Since $u$ is continuous and $F$ is upper semicontinuous, for any $z\in\Omega$ we have
\begin{equation}\label{usc}
    u(z)=\lim_{z_i\rightarrow z}u(z_i)\leq \limsup_{z_i\rightarrow z} F(z_i) \leq F(z).
\end{equation}
This proves (\ref{claim1}). In fact we proved 
\begin{equation*}
    \lim_{n\rightarrow\infty} \frac{1}{n}\log |f_n(z_i)| = u(z_i)
\end{equation*}
for all $z_i$. \par

\emph{Step 5.---} Lastly we prove (\ref{claim2}). For convenience we collect what we have proved so far.
\begin{align*}
&(1)\; \frac{1}{n}\log |f_n|\; \text{is locally uniformly bounded from above ;}\\
&(2)\;\limsup_{n\rightarrow\infty}\frac{1}{n}\log |f_n(z)|\leq u(z) \quad\forall z\in\Omega \,;\\
&(3)\;\lim_{n\rightarrow\infty} \frac{1}{n}\log |f_n(z_i)| = u(z_i) \quad\forall z_i.
\end{align*}
 By (1) and (3), the sequence $\frac{1}{n}\log |f_n(z)|$ is compact in $L^1_{loc}(\Omega)$. By (2) and (3), every subsequential $L^1_{loc}(\Omega)$ limit is $\leq u$ with equality on a dense subset $(z_i)_{i=1}^{\infty}$. By the same argument as in (\ref{usc}), all limits are equal to $u$. Hence $\frac{1}{n}\log|f_n(z)|$ converges to $u$ in $L^1_{loc}(\Omega)$. This argument is due to H\"{o}rmander \cite[Proof of Theorem 4.2.13]{hormander2007notions}.
\end{proof}

As a corollary we prove a well-known fact that any psh function $u$ on $\Omega$ can be approximated by functions of the form $N^{-1}\log|f|$ with $N\in\mathbb{N}$ and $f\in\mathcal{O}(\Omega)$.
\begin{corollary}[\protect{\cite[Theorem 4.2.13]{hormander2007notions}}] \label{horcor}
    Let $\Omega\subset\mathbb{C}^N$ be a bounded pseudoconvex domain. Then $\{m^{-1}\log|f| : f\in\mathcal{O}(\Omega), m\in\mathbb{N}\}$ is dense in $PSH(\Omega)$ equipped with the $L^1_{loc}(\Omega)$ topology.
\end{corollary}
\begin{proof}
We need to show that for arbitrary $u\in PSH(\Omega)$, there exists a sequence of $h_n\in\mathcal{O}(\Omega)$ and $m_n\in\mathbb{N}$ such that $m_n^{-1}\log|h_n|$ converges to $u$ in $L^1_{loc}(\Omega)$. Since continuous psh functions are dense in $PSH(\Omega)$ in $L^1_{loc}(\Omega)$-topology, we can assume $u$ is continuous without loss of generality. Then Theorem \ref{main} provides such a sequence almost surely.
\end{proof}
When $N=1$ and $\Omega$ is simply connected, one can give a more elementary proof by approximating the measure $\frac{m}{2\pi}\Delta u$ by integer point masses. See the remark after Theorem 4.2.13 in \cite{hormander2007notions}.

Now we estimate the tail probability of random variables induced by classical linear functionals on $PSH(\Omega)$ such as $M_v(z,r):=\fint_{B(z,r)}v \,d\lambda=\lambda(B(z,r))^{-1}\int_{B(z,r)}v \,d\lambda$, $C_v(z,r):=\fint_{\partial B(z,r)}v \,d\sigma=\sigma(\partial B(z,r))^{-1}\int_{\partial B(z,r)}v \,d\sigma$ and $v_{\varepsilon}(z):= v*\rho_{\varepsilon}(z)$ where $\rho_{\varepsilon}$ is a standard mollifier. For simplicity we only prove the following theorem with $M_v(z,r)$, but it is clear that the same argument can be applied to other functionals. It is a straightforward consequence of Proposition \ref{dem}, subaveraging, and the proof of \emph{step 1} of Theorem \ref{main}.

\begin{proposition}\label{subexpthm}
    Let $u$ and $(f_n)_n$ be as in Theorem \ref{main}. Given $z\in\Omega$ and $\varepsilon>0$, there exists a constant $0<R< d(z,\partial\Omega)$ such that for all $0<r<R$,
    \begin{equation}\label{subexp}
        \mathbb{P}\left( \left| \fint_{B(z,r)}\frac{1}{n}\log |f_n|\,d\lambda-\fint_{B(z,r)}u\,d\lambda \right| >\varepsilon \right)\leq C e^{-Dn}
    \end{equation}
     holds for all $n\in\mathbb{N}$ with some constants $C,D>0$.
\end{proposition}
\begin{proof}
We only need to prove (\ref{subexp}) for large enough $n$. By Proposition \ref{dem} we can let $n$ large enough so that
\begin{equation*}
    \left| \fint_{B(z,r)}\frac{1}{n}\log S_n \,d\lambda-\fint_{B(z,r)}u\,d\lambda \right| < \frac{\varepsilon}{2}
\end{equation*}
for a given $r>0$ which will be specified later. Then
\begin{equation}\label{summand}
\begin{gathered}
    \mathbb{P}\left( \left| \fint_{B(z,r)}\frac{1}{n}\log |f_n|-u \, d\lambda \right| >\varepsilon \right)
    \leq \mathbb{P}\left( \left| \fint_{B(z,r)}\frac{1}{n}\log |f_n|-\frac{1}{n}\log S_n \, d\lambda \right| >\frac{\varepsilon}{2} \right) \\
    \leq \mathbb{P}\left(  \fint_{B(z,r)}\frac{1}{n}\log \frac{|f_n|}{S_n} \, d\lambda >\frac{\varepsilon}{2} \right)
+\,\mathbb{P}\left(  \fint_{B(z,r)}\frac{1}{n}\log \frac{|f_n|}{S_n} \, d\lambda < -\frac{\varepsilon}{2} \right).
\end{gathered}
\end{equation}

We estimate each summand in (\ref{summand}) separately. For the first summand, (\ref{step1}) shows that there exists $0<R_1<d(z,\partial\Omega)$ such that for $0<r<R_1$,
\begin{equation*}
\begin{gathered}
\mathbb{P}\left(  \fint_{B(z,r)}\frac{1}{n}\log \frac{|f_n|}{S_n} \, d\lambda >\frac{\varepsilon}{2} \right)\leq\mathbb{P}\left(  \sup_{B(z,r)}\frac{1}{n}\log \frac{|f_n|}{S_n} >\frac{\varepsilon}{2} \right) \\
=\mathbb{P}\left(  \sup_{B(z,r)} \frac{|f_n|}{S_n} >e^{{\varepsilon n}/2} \right) \leq C_1 e^{-D_1 n}
\end{gathered}
\end{equation*}
holds with some constants $C_1, D_1 >0$. 

Now we move on to the second summand. Let $\delta>0$ and $0<R_2<d(z,\partial\Omega)$ so that $\varepsilon':=\varepsilon/2-\delta>0$ and 
\begin{equation*}
    \left| u(z) - \fint_{B(z,r)}u\,d\lambda \right|<\frac{\varepsilon'}{3}
\end{equation*}
for $0<r<R_2$.
Since $\frac{1}{n}\log S_n$ converges to $u$ pointwisely and in $L^1_{loc}(\Omega)$ by Proposition \ref{dem}, for each $0<r<R_2$ we can let $n$ large enough so that 
\begin{equation*}
    \left| \frac{1}{n}\log S_n(z)-u(z) \right|, \left| \fint_{B(z,r)}\frac{1}{n}\log S_n(z')-u(z')\,d\lambda \right| < \frac{\varepsilon'}{3}
\end{equation*}
and hence 
\begin{equation}\label{varepsilon}
    \left|\frac{1}{n}\log S_n(z) - \fint_{B(z,r)}\frac{1}{n}\log S_n\,d\lambda \right|<\varepsilon'.
\end{equation}

Let $\mathcal{A}$ and $\mathcal{B}$ be the events defined by 
\begin{align*}
\mathcal{A}&:=\{\frac{1}{n}\log |f_n(z)|-\fint_{B(z,r)}\frac{1}{n}\log S_n \, d\lambda < -\frac{\varepsilon}{2}\} \,; \\
\mathcal{B}&:= \{\frac{1}{n}\log |f_n(z)|>\frac{1}{n}\log S_n(z)-\delta \}.
\end{align*}
Then by the submean-value inequality, we have
\begin{equation*}
\mathbb{P}\left( \fint_{B(z,r)}\frac{1}{n}\log |f_n|-\frac{1}{n}\log S_n \, d\lambda < -\frac{\varepsilon}{2} \right) \leq \mathbb{P}(\mathcal{A})\leq
\mathbb{P}(\mathcal{A}\cap\mathcal{B})+\mathbb{P}(\mathcal{A}\cap\mathcal{B}^c).
\end{equation*}
For large enough $n$, (\ref{varepsilon}) then implies
\begin{equation*}
    \mathbb{P}(\mathcal{A}\cap\mathcal{B})\leq \mathbb{P}\left(\frac{1}{n}\log S_n(z) - \fint_{B(z,r)}\frac{1}{n}\log S_n\,d\lambda <\delta-\frac{\varepsilon}{2}\right)=0.
\end{equation*}
$\mathbb{P}(\mathcal{B}^c)$ can be bounded by an argument similar to that in Lemma \ref{inf}:
\begin{equation*}
\mathbb{P}(\mathcal{B}^c)=\mathbb{P}\left(\frac{|f_n(z)|}{S_n(z)}\leq e^{-\delta n}\right)\leq C_2 e^{-D_2 n}
\end{equation*}
for some constants $C_2, D_2>0$. This completes the proof.
\end{proof}
\begin{remark}
    In fact, by \eqref{localuniform} the dependence of $C, D$ on $0<r<R$ only comes from finitely many small $n$ which may differ for each $r$. Hence Proposition \ref{subexpthm} holds asymptotically with the same $C, D$ for all $0<r<R$ as $n\rightarrow\infty$.
\end{remark}

\section{Convergence of zeros}

For a holomorphic map $F=(f^1,...,f^q):\Omega\to\mathbb{C}^q$ that is transverse to the origin, we write
\begin{equation*}
  Z_{f^1,...,f^q}:=[F^{-1}(0)]
\end{equation*}
for the current of integration over the resulting smooth submanifold. We first justify the almost sure transversality needed for applications in zeros of GAFs. 

\begin{lemma}\label{transversality}
Fix $n\geq1$, let $1\leq k\leq N$, and let $f_n^1,...,f_n^k$ be independent copies of the GAF in Definition \ref{GAF}.  With probability one, for every nonempty subset $I\subset\{1,...,k\}$ the map
\begin{equation*}
 F_{n,I}:=(f_n^j)_{j\in I}:\Omega\longrightarrow\mathbb{C}^{|I|}
\end{equation*}
is transverse to $0$. Consequently, $F_{n,I}^{-1}(0)$ is either empty or a smooth complex submanifold of codimension $|I|$.
\end{lemma}

\begin{proof}
Fix $p\in\Omega$.  By \eqref{nonzero}, $B_n(p,p)=S_n(p)^2>0$. Hence there is a unit vector $h\in\mathcal{H}(nu)$ such that $h(p)\neq0$. Choose an open neighborhood $U\subset\subset\Omega$ on which $h$ has no zeros. Complete $h$ to an orthonormal basis of $\mathcal{H}(nu)$.  The law of the GAF is invariant under a unitary change of orthonormal basis; hence, in the corresponding Gaussian expansions, for every $j\in I$ we may write
\begin{equation*}
 f_n^j=a_jh+g_j,
\end{equation*}
where $(a_j)_{j\in I}$ are independent standard complex Gaussian variables and are independent of the random holomorphic functions $(g_j)_{j\in I}$.

Condition on a realization $g=(g_j)_{j\in I}$ and put $r=|I|$. Consider the smooth map
\begin{equation*}
 \Phi_g : U\times\mathbb{C}^r \longrightarrow\mathbb{C}^r, \qquad \Phi_g(z,a)=\left( g_j(z)+a_jh(z) \right)_{j\in I}.
\end{equation*}
Since
\begin{equation*}
 D_a \Phi_g(z,a)=h(z) \operatorname{Id}_{\mathbb{C}^r},
\end{equation*}
$\Phi_g$ is a submersion. Thus $M_g:=\Phi_g^{-1}(0)$ is a smooth manifold. Let $\pi_g:M_g\to\mathbb{C}^r$ be the restriction of the projection onto the
parameter factor.  At $(z,a)\in M_g$,
\begin{equation*}
 T_{(z,a)}M_g =\{(\xi,b):D_z\Phi_g(z,a)\xi+h(z)b=0\}.
\end{equation*}
It follows that $D\pi_g$ is surjective at $(z,a)$ if and only if $D_z\Phi_g(z,a):T_zU\to\mathbb{C}^r$ is surjective.  Equivalently, $a$ is a regular value of $\pi_g$ if and only if the map $z\mapsto\Phi_g(z,a)$ is transverse to $0$ on $U$. By Sard's theorem, for every fixed $g$, the set of exceptional parameters $a\in\mathbb{C}^r$ has Lebesgue measure zero, and hence Gaussian measure zero. By applying this argument to a compact exhaustion of $\Omega$ and then using Fubini's theorem, we obtain the lemma.

\end{proof}

We now prove the almost sure convergence of the simultaneous zero currents of GAFs. In fact, the proof below does not use Theorem \ref{main} and is logically independent of it.

\begin{theorem}\label{zero}
For every $n\geq1$, let $f_n^1,...,f_n^N$ be independent copies of the GAF in Definition \ref{GAF}. Then for every $1\leq k\leq N$,
\begin{equation}\label{expzero}
\frac{1}{n^k}\mathbb{E} Z_{f_n^1,...,f_n^k}=\left( dd^c \frac{1}{n}\log S_n \right)^k \longrightarrow (dd^cu)^k.
\end{equation}
Moreover,
\begin{equation}\label{justzero}
\frac1{n^k}Z_{f_n^1,...,f_n^k}\longrightarrow (dd^cu)^k
\end{equation}
weakly as currents on $\Omega$, almost surely. 
\end{theorem}

\begin{proof}
Since $u$ is continuous, from Demailly's estimate \eqref{dem} it is easy to see that
\begin{equation}\label{localuniform}
 u_n\longrightarrow u \quad \text{locally uniformly on}\;\Omega.
\end{equation}

Using Lemma \ref{transversality} for all $n$ and all subsets of $\{1,...,N\}$, define
\begin{equation*}
 R_{n,q}:=\frac1{n^q}Z_{f_n^1,...,f_n^q}, \qquad R_{n,0}:=1.
\end{equation*}
The currents $R_{n,q}$ are well defined with probability one, and the corresponding zero sets are smooth submanifolds. Finally, let
\begin{equation*}
 v_{n,q}:=\frac{1}{n}\log|f_n^q|, \qquad x_{n,q}:=v_{n,q}-u_n =\frac{1}{n}\log\frac{|f_n^q|}{S_n}.
\end{equation*}

We proceed by induction. Let $Y_{n,q-1}=\{f_n^1=\cdots=f_n^{q-1}=0\}$.  On each nonempty component of $Y_{n,q-1}$, the restriction $f_n^q|_{Y_{n,q-1}}$ is not identically zero; otherwise $df_n^q$ would lie in the span of
$df_n^1,...,df_n^{q-1}$, hence contradicting the transversality of
$(f_n^1,...,f_n^q)$. By Poincar\'e--Lelong formula on $Y_{n,q-1}$, we have
\begin{equation*}
    Z(f_n^q|_{Y_{n,q-1}})=dd^c \left(\log |f_n^q|_{Y_{n,q-1}}|\right).
\end{equation*}
Pushforward by inclusion $\iota: Y_{n,q-1} \hookrightarrow \Omega$ gives
\begin{equation*}
 dd^c\left(\log|f_n^q|[Y_{n,q-1}]\right)=\iota_* dd^cf_n^q=\iota_* Z_{Y_{n,q-1}}(f_n^q)= [Y_{n,q-1}\cap\{f_n^q=0\}].
\end{equation*}
Note that the multiplication of currents on the left hand side is well defined. Dividing by $n^q$ yields
\begin{equation}\label{upluserr}
 R_{n,q}=dd^c\left( v_{n,q}R_{n,q-1} \right) =dd^c\left( (u_n+x_{n,q})R_{n,q-1} \right).
\end{equation}
By \eqref{localuniform}, $u_n$ is well-behaved. Hence we need to control the $x_{n,q}$ term in the proofs of both \eqref{expzero} and \eqref{justzero}.

First we prove \eqref{expzero}. Let $\xi$ be a standard complex Gaussian variable. For every fixed $z$,
\begin{equation*}
 \frac{f_n^q(z)}{S_n(z)}\sim\xi.
\end{equation*}
Consequently, for finite constants
$c_0:=\mathbb{E}\log|\xi|$ and $C_0:=\mathbb{E}|\log|\xi||$,
\begin{equation*}
 \mathbb{E} x_{n,q}(z)=\frac{c_0}{n}, \qquad \mathbb{E}|x_{n,q}(z)|=\frac{C_0}{n}.
\end{equation*}
Since $f_n^q$ is independent of $R_{n,q-1}$, the expectation is easy to obtain:
\begin{equation*}
 \mathbb{E}_{f_n^q}\left[x_{n,q}R_{n,q-1}\right] =\frac{c_0}{n}R_{n,q-1}
\end{equation*}
as currents. Taking conditional expectation in \eqref{upluserr} using closedness of $R_{n,q-1}$, and then inducting on $q$, we obtain the exact identity
\begin{equation}\label{exact-mean-current}
 \mathbb{E} R_{n,q}=(dd^c u_n)^q, \quad 0\leq q \leq N.
\end{equation}

By \eqref{localuniform} and continuity of the Bedford--Taylor product \cite{BT82} under locally uniform convergence of locally bounded psh potentials (see \cite{guedj2017degenerate}), we have
\begin{equation}\label{BT}
 (dd^c u_n)^q\longrightarrow (dd^c u)^q, \quad 0\leq q\leq N.
\end{equation}
This proves \eqref{expzero}. In particular, if $\beta:=dd^c|z|^2$, $\chi\in C_c^\infty(\Omega)$, and $\chi\geq0$, then
\begin{equation}\label{mass}
 \sup_n\int_\Omega\chi (dd^c u_n)^q\wedge\beta^{N-q}<\infty.
\end{equation}

Next we prove \eqref{justzero} by induction on $q$. Recall that we need to prove $R_{n,q}\to (dd^c u)^q$ almost surely, beginning with $R_{n,0}=1$. Suppose that
\begin{equation}\label{induction}
 R_{n,q-1}\longrightarrow (dd^c u)^{q-1}
\end{equation}
almost surely.  Fix $K\subset\subset\Omega$ and choose $\chi\in C_c^\infty(\Omega)$ with $\chi\geq0$ and $\chi=1$ on a neighborhood of $K$. To control $x_{n,q}$, define the nonnegative random variable
\begin{equation*}
 Y_{n,\eta}:= \int_\Omega \chi |x_{n,q}| \mathbf{1}_{\{|x_{n,q}|>\eta\}} R_{n,q-1}\wedge\beta^{N-q+1}
\end{equation*}
for $\eta>0$. Set
\begin{equation*}
 b_n(\eta):= \mathbb{E}\left[ \left|\frac{1}{n}\log|\xi|\right| \mathbf{1}_{\{|n^{-1}\log|\xi||>\eta\}} \right].
\end{equation*}
Using the radial density $2re^{-r^2}\,dr$ of $|\xi|$, the contribution from
$0<r<e^{-\eta n}$ is bounded by
\begin{equation*}
 \frac2n\int_0^{e^{-\eta n}}(-\log r)r\,dr =e^{-2\eta n}\left(\eta+\frac1{2n}\right).
\end{equation*}
For $r>e^{\eta n}$, the inequality $\log r\leq r^2/2$ gives
\begin{equation*}
 \frac2n\int_{e^{\eta n}}^\infty(\log r)re^{-r^2}\,dr \leq\frac{e^{2\eta n}+1}{2n}e^{-e^{2\eta n}}.
\end{equation*}
Therefore
\begin{equation}\label{logtail}
 \sum_{n=1}^\infty b_n(\eta)<\infty \quad\text{for every}\;\eta>0.
\end{equation}

By \eqref{exact-mean-current} and \eqref{mass}, we have
\begin{align*}
 \mathbb{E}Y_{n,\eta}
 &=b_n(\eta)\int_\Omega\chi\,\mathbb{E} R_{n,q-1}\wedge\beta^{N-q+1}\\
 &=b_n(\eta)\int_\Omega\chi dd^c u_n^{q-1}\wedge\beta^{N-q+1} \leq C_{K,q}b_n(\eta).
\end{align*}
It follows from \eqref{logtail} and Tonelli's theorem that $\sum_nY_{n,\eta}<\infty$ almost surely, and hence
\begin{equation*}
 Y_{n,\eta}\longrightarrow0 \quad\text{almost surely}.
\end{equation*}
On the probability one event on which \eqref{induction} holds, positivity and weak convergence give
\begin{equation}\label{Mkq}
 M_{K,q-1}:= \sup_n\int_\Omega\chi R_{n,q-1}\wedge\beta^{N-q+1}<\infty.
\end{equation}
Consequently,
\begin{equation*}
 \int_K|x_{n,q}|R_{n,q-1}\wedge\beta^{N-q+1} \leq \eta M_{K,q-1}+Y_{n,\eta}.
\end{equation*}
Apply this for $\eta=1/m$, intersect the corresponding countably many probability-one events, and let $m\to\infty$. We obtain
\begin{equation*}
 \int_K|x_{n,q}|R_{n,q-1}\wedge\beta^{N-q+1} \longrightarrow 0.
\end{equation*}
Positivity of $R_{n,q-1}$ implies that the trace measure dominates the total variation \cite[Corollary 2.24]{guedj2017degenerate}. Hence we prove
\begin{equation}\label{xR}
 x_{n,q}R_{n,q-1}\longrightarrow 0
\end{equation}
as currents.

Furthermore, \eqref{localuniform}, \eqref{induction}, and \eqref{Mkq} imply
\begin{equation}\label{uR}
 u_nR_{n,q-1}- u(dd^c u)^{q-1}=(u_n-u)R_{n,q-1}+u(R_{n,q-1}-(dd^c u)^{q-1})\longrightarrow 0.
\end{equation}
For the second term, we need to uniformly approximate $u$ by smooth functions.

Finally, \eqref{upluserr}, \eqref{xR}, and
\eqref{uR} yield
\begin{equation*}
 R_{n,q}=dd^c \left( (u_n+x_{n,q})R_{n,q-1} \right) \longrightarrow dd^c\left(u(dd^c u)^{q-1}\right)=(dd^c u)^q.
\end{equation*}
This completes the induction step.
\end{proof}

In complex dimension $N=1$, the statement says that the normalized point process $n^{-1}Z_{f_n}$ converges to $dd^cu=\frac1{2\pi}\Delta u$ in the sense of distributions, consistently with Theorem~1 of \cite{Sodin}.

\bibliographystyle{abbrv}
\bibliography{References}

\end{document}